\documentclass[11pt,openbib]{article} 
\usepackage[T1]{fontenc}
\usepackage{amsfonts,amsthm}
\usepackage{amssymb,amsmath}
\usepackage{epic}
\newcommand{\R}{\mathbb{R}}

\newcommand{\dee}{\mathrm{d}}

\newcommand{\vvee}{\mbox{\tiny $\vee $}}
\newcommand{\comp}{\raisebox{0pt}{$\scriptstyle\circ \, $}}

\newcommand{\setrule}{\, \mathop{\rule[-4pt]{.5pt}{13pt}\, }\nolimits}

\newcommand{\bigspace}{\bigskip\par\noindent}

\allowdisplaybreaks
\begin{document}
\mbox{}\vspace{.1in}
\begin{center}
{\bf \Large Differential spaces in integrable Hamiltonian systems}
\vspace{.05in}\par\noindent Richard Cushman and J\k{e}drzej \'{S}niatycki\footnotemark 
\end{center} 
\footnotetext{\parbox[t]{4in}{email: rcushman@gmail.com and sniatycki@gmail.com \\ 
Department of Mathematics and Statistics, University of Calgary} } \bigskip 

In this paper we use differential space to obtain some new results for completely integrable Hamiltonian systems. Differential spaces were introduced by Sikorski \cite{sikorski67} in 1967 with a comprehensive treatment given in his book \cite{sikorski72} (in Polish). The theory of differential spaces allows the differential geometric study of arbitrary subsets of ${\R }^n$ \cite{sniatycki13}. In particular it allows one to study the singularities of classical integrable systems \cite{cushman-bates}. In a recent paper \cite{seppe-vungoc} on completely integrable systems 
the authors avoid the theory of differential spaces using ad hoc definitions. They obtain weaker results than 
the ones given here. \bigskip 

A \emph{completely integrable Hamiltonian system} is a triple 
$(M, \omega, \mathbf{F} =(f_1, \ldots ,$ $f_n) \big)$, 
where 1) $(M, \omega )$ is a $2n$-dimensional smooth symplectic manifold; 2) the Hamiltonian vector fields 
$X_{f_i}$ for $1 \le i \le n$ are complete; 3) the components of the \emph{integral mapping} 
\begin{displaymath}
\mathbf{F}:M \rightarrow \R : x \mapsto \big( f_1(x), \ldots , f_n(x) \big) 
\end{displaymath} 
are in pairwise involution, that is, $L_{X_{f_j}} f_i = \{ f_i, f_j \} =0$; and 4) the components are 
functionally independent almost everywhere, that is, the rank of $D\mathbf{F}(x)$ is $n$ for 
all points $x$ in a dense open subset $U$ of $M$. \medskip 

Let $M/X_{\mathbf{F}}$ be the space of orbits of the family of Hamiltonian vector fields 
$X_{\mathbf{F}} = {\{ X_{f_i} \} }^n_{i=1}$ on $(M,\omega )$, see \cite[Chpt 3, p.44]{sniatycki13}. \medskip

\noindent \textbf{Claim 1} Each orbit of a family of Hamiltonian vector fields $X_{\mathbf{F}}$ 
associated to the completely integrable system $(M, \omega , \mathbf{F})$ is an immersed submanifold of 
$M$. \medskip 

\noindent \textbf{Proof.} This follows from the first Sussmann theorem \cite[theorem 4.1 p.179]{sussmann}. See also  
\'{S}niatycki \cite[theorem 3.4.5 p.46]{sniatycki13}.  \hfill $\square $ \medskip 

The \emph{commutant} of the integrable system $(M, \omega, \mathbf{F})$ is 
the set ${C^{\infty}(M)}^{\mathbf{F}}$ of smooth functions on $M$, which Poisson commute with $f_i$ for 
$1\le i \le n$, that is, $g \in {C^{\infty}(M)}^{\mathbf{F}}$ if and only if 
$\{ g, f_i \} = 0$ for every $1 \le i \le n$. The notion of commutant was introduced in 
\cite{seppe-vungoc}. Equivalently, \medskip  

\noindent \textbf{Claim 2} The commutant ${C^{\infty}(M)}^{\mathbf{F}}$ is a differential structure 
on the space of orbits $M/X_{\mathbf{F}}$.  \medskip 

\noindent \textbf{Proof.} Let $\pi : M \rightarrow M/X_{\mathbf{F}}$ be the canonical projection map. By 
definition $h \in C^{\infty}(M/X_{\mathbf{F}})$ if and only if ${\pi }^{\ast }h \in {C^{\infty}(M)}^{\mathbf{F}}$. 
We show that $C^{\infty}(M/X_{\mathbf{F}})$ satisfies the conditions for a differential structure 
\cite[Chpt 2, p.15]{sniatycki13}. \medskip 

Consider the topology on $M/X_{\mathbf{F}}$ generated by the subbasis 
\begin{displaymath}
\big\{ h^{-1}(I) \setrule \mbox{$I$ and open interval in $\R $ and $h \in
C^{\infty}(M/X_{\mathbf{F}})$}\big\} . 
\end{displaymath}
This topology satisfies condition $1$. Condition 2 is automatic. 
To verify condition 3 let $g: M/X_{\mathbf{F}} \rightarrow \R $ 
be a function such that for every $\overline{x} \in M/X_{\mathbf{F}}$ there are functions 
$f, f_1, \ldots , f_n  \in C^{\infty}(M/X_{\mathbf{F}})$ and open intervals 
$I_1, \ldots , I_{\ell } \subseteq \R $ such that 
\begin{equation}
\overline{x} \in U_{\overline{x}} = f^{-1}_1(I_1) \cap \cdots f^{-1}_{\ell }(I_{\ell }) 
\label{eq-one}
\end{equation}
and 
\begin{equation} 
g|_{U_{\overline{x}}} = f|_{U_{\overline{x}}}.
\label{eq-two}
\end{equation}
Since $f_1, \ldots , f_{\ell } \in C^{\infty}(M/X_{\mathbf{F}})$, it follows that 
\begin{align} 
{\pi }^{-1}(U_{\overline{x}}) & = {\pi }^{-1}\big( f^{-1}_1(I_1) \cap \cdots \cap f^{-1}_{\ell }(I_{\ell }) \big)
\notag \\
& = {\pi }^{-1}\big( f^{-1}_1(I_1)\big) \cap \cdots \cap {\pi }^{-1} \big( f^{-1}_{\ell }(I_{\ell }) \big) 
\notag \\
& = (f_1 \comp \pi )^{-1}(I_1) \cap \cdots \cap (f_{\ell } \comp \pi )^{-1}(I_{\ell }) \notag 
\end{align}
is open in $M$. Hence the topology of the orbit space $M/X_{\mathbf{F}}$ is coarser than that of $M$. 
Equation (\ref{eq-two}) implies that $\{ {\pi }^{-1}(U_{\overline{x}}) \setrule \overline{x} \in M/X_{\mathbf{F}} \} $ 
is an open covering of $M$ such that ${\pi }^{\ast }\big( g|_{U_{\overline{x}}}\big) $ is smooth and is 
preserved by the vector fields $X_{f_1}, \ldots X_{f_n}$ restricted to ${\pi }^{-1}(U_{\overline{x}})$. Hence  
${\pi }^{\ast }g \in {C^{\infty}(M)}^{\mathbf{F}}$. \hfill $\square $ \medskip 

\noindent \textbf{Corollary 2a} The orbit space $M/X_{\mathbf{F}}$ with its differential structure 
${C^{\infty}(M)}^{\mathbf{F}}$ is the differential space $\big( M/X_{\mathbf{F}}, C^{\infty}(M/X_{\mathbf{F}})\big)$. \medskip 

Because the hypothesis that the vector fields in the family $X_{\mathbf{F}}$ are complete is not 
used in the proofs above, we introduce a weaker notion of \emph{integrable Hamiltonian system}, 
which satsifies hypotheses 1), 3) and 4) of the definition of completely integrable Hamiltonian system. 
We have proved \medskip 

\noindent \textbf{Corollary 2b} Let $\big( M, \omega , \mathbf{F} = (f_1, \ldots , f_n) \big) $ be an 
integrable Hamiltonian system. Then each orbit of $X_{\mathbf{F}}$ is an immersed submanifold of 
$M$. The orbit space $M/X_{\mathbf{F}}$ with its differential structure $C^{\infty}(M/X_{\mathbf{F}})$ is a differential 
space and the projection mapping $\pi :M \rightarrow M/X_{\mathbf{F}}$ is smooth. \medskip 

Let $(M, \omega , \mathbf{F})$ and $(M, \omega , \mathbf{G})$ be two integrable Hamiltonian systems 
on $(M, \omega )$. Following Seppe and Vu Ngoc \cite{seppe-vungoc}, we say that these integrable systems are equivalent, that is, $\mathbf{F} \sim \mathbf{G}$, if and only if ${C^{\infty}(M)}^{\mathbf{F}} = {C^{\infty}(M)}^{\mathbf{G}}$. Let $X_{{C^{\infty}(M)}^{\mathbf{F}}} $ be 
the set of all Hamiltonian vector fields $X_f$ on $(M, \omega )$ where $f$ lies in the commutant 
${C^{\infty}(M)}^{\mathbf{F}}$. Observe that the orbit of the families $X_{\mathbf{F}}$ and 
$X_{{C^{\infty}(M)}^{\mathbf{F}}}$ through each point $x$ in $M$ coincide. Hence, 
$\mathbf{F} \sim \mathbf{G}$ is equivalent to saying that on $M$ the orbits of 
the families $X_{\mathbf{F}}$ and $X_{\mathbf{G}}$ coincide. \medskip 

As in \cite{seppe-vungoc}, we say that two integrable systems $(\mathbf{F}, M, \omega )$ and 
$({\mathbf{F}}', M', {\omega }' )$ 
are \emph{symplectically equivalent} if and only if there is a symplectic diffeomorphism 
$\varphi : M \rightarrow M'$ such that $\mathbf{F} = {\varphi }^{\ast }{\mathbf{F}}'$. In other words, 
${C^{\infty}(M)}^{\mathbf{F}} = {C^{\infty}(M)}^{{\varphi }^{\ast }{\mathbf{F}}'}$. Because the diffeomorphism $\varphi $ induces the map ${\varphi }^{\ast }: C^{\infty}(M') \rightarrow C^{\infty}(M)$ being symplectically equivalent means that we have ${C^{\infty}(M)}^{\mathbf{F}} = $ 
${\varphi }^{\ast } \big( {C^{\infty}(M')}^{{\mathbf{F}}'} \big) $. If $(\mathbf{F}, M, \omega )$ and
$({\mathbf{F}}', M', {\omega }' )$ are symplectically equivalent by the diffeomorphism 
$\varphi : M \rightarrow M'$, then $\varphi $ 
induces the diffeomorphism of differential spaces
\begin{displaymath}
{\varphi }^{\vvee}: \big( M/X_{\mathbf{F}}, C^{\infty}(M/X_{\mathbf{F}}) \big) \rightarrow 
\big( M'/{\mathbf{F}}', C^{\infty}(M'/{\mathbf{F}}') \big) 
\end{displaymath}
defined by sending the orbit of $X_{\mathbf{F}}$ through $x \in M$ to the orbit of $X_{{\mathbf{F}}'}$ in 
$M'$ through $\varphi (x)$. Moreover the following diagram commutes. 
\bigspace
\setlength{\unitlength}{.5mm}
\hspace{1.25in}\begin{tabular}{l}
\begin{picture}(90,55)
\put(28,51){\makebox(0,0)[tl]{$\mathrm{M}$}}
\put(37,48){\vector(1,0){31}}
\put(70,51.5){\makebox(0,0)[tl]{$M'$}}
\put(48,49){\makebox(0,0)[bl]{$\varphi $}}
\put(31,44){\vector(0,-1){18}}
\put(33,34){\makebox(0,0)[bl]{$\pi $}}
\put(73,44){\vector(0,-1){19}}
\put(75,34){\makebox(0,0)[bl]{${\pi }'$}}
\put(46,20){\vector(1,0){20}}
\put(22,23){\makebox(0,0)[tl]{$M/X_{\mathbf{F}}$}}
\put(53,21){\makebox(0,0)[bl]{${\varphi }^{\vvee}$}}
\put(67,23){\makebox(0,0)[tl]{$M/{\mathbf{F}}'$}}
\end{picture} \\
\end{tabular} \vspace{-10pt}
  
\noindent \textbf{Claim 3} If $(\mathbf{F}, M, \omega )$ and $({\mathbf{F}}', M', {\omega }' )$ are symplectically equivalent by the symplectic diffeomorphism $\varphi $ and $(M', {\omega }', {\mathbf{F}}' )$ is completely integrable, then $(M, \omega , \mathbf{F} )$ is completely integrable. \medskip 

\noindent \textbf{Proof.} Suppose that $\varphi : M \rightarrow M'$ is a symplectic diffeomorphism 
such that $\mathbf{F} = {\varphi }^{\ast }{\mathbf{F}}'$. Then 
\begin{displaymath}
X_{\mathbf{F}} = (X_{f_1}, \ldots , \, X_{f_n} ) = (X_{{\varphi }^{\ast }f'_1}, \ldots , 
X_{{\varphi }^{\ast }f'_n} ) = X_{{\varphi }^{\ast }{\mathbf{F}}'}. 
\end{displaymath}
So if the Hamiltonian vector fields in $X_{{\varphi }^{\ast }{\mathbf{F}}'}$ are complete then those in 
$X_{\mathbf{F}}$ are also. \hfill $\square $ \medskip 

We now look at the relation between the space of orbits $M/X_{\mathbf{F}}$ of the family $X_{\mathbf{F}}$ of 
Hamiltonian vector fields associated to the integrable Hamiltonian system $(M, \omega ,\mathbf{F} )$ and 
the image $\mathbf{F}(M)$ of its integral map. Since $\mathbf{F}(M) \subseteq {\R }^n$, we have 
a differential structure $C^{\infty}_i(\mathbf{F}(M))$ on $\mathbf{F}(M)$ defined by $g \in 
C^{\infty}_i(\mathbf{F}(M))$ if and only if for every $y \in \mathbf{F}(M)$ there is an open neighborhood 
$V_y$ of $y$ in ${\R }^n$ and a function $G_y \in C^{\infty}({\R }^n)$ such that 
$g|_{V_y \cap \mathbf{F}(M)} = G_y|_{V_y\cap \mathbf{F}(M)}$. The differential space 
$\big( \mathbf{F}(M), C^{\infty}_i(\mathbf{F}(M)) \big) $ 
is a differential subspace of $\big( {\R }^n, C^{\infty}({\R }^n) \big) $. \medskip 

\noindent \textbf{Claim 4} The integral map 
\begin{displaymath}
\mathbf{F}: \big( M, C^{\infty}(M) \big) \rightarrow \big( \mathbf{F}(M), C^{\infty}_i(\mathbf{F}(M)) \big) 
\end{displaymath}
is a smooth mapping of differential spaces. \medskip 

\noindent \textbf{Proof.} Suppose that $g \in C^{\infty}_i(\mathbf{F}(M))$. We need to show that 
${\mathbf{F}}^{\ast }g \in C^{\infty}(M)$. By definition for every $y \in \mathbf{F}(M)$ there is an open 
neighborhood $V_y$ of $y$ in ${\R }^n$ and a function $G_y \in C^{\infty}({\bf \R }^n)$ such that 
$g|_{V_y \cap \mathbf{F}(M)} = G_y|_{V_y \cap \mathbf{F}(M)}$. Hence 
\begin{equation}
({\mathbf{F}}^{\ast }g)|_{{\mathbf{F}}^{-1}(V_y\cap \mathbf{F}(M))} = 
({\mathbf{F}}^{\ast }G_y)|_{{\mathbf{F}}^{-1}(V_y\cap \mathbf{F}(M))}. 
\label{eq-three}
\end{equation}
Since $G_y \in C^{\infty}({\R }^n)$ and $\mathbf{F}: M \rightarrow \mathbf{F}(M) \subseteq {\R }^n$ is 
smooth, it follows that ${\mathbf{F}}^{\ast }G_y = G_y \comp \mathbf{F}$ is smooth. Thus for every 
$x \in M$ there is a $y= \mathbf{F}(x) \in \mathbf{F}(M)$, an open neighborhood 
${\mathbf{F}}^{-1}(V_y)$ of $x$ in $M$, and a function ${\mathbf{F}}^{\ast }G_y \in C^{\infty}(M)$ such that 
equation (\ref{eq-three}) holds. Hence ${\mathbf{F}}^{\ast }g \in C^{\infty}(M)$. \hfill $\square $ \medskip 

For each $x \in M$ let $L_x$ be the connected component of the fiber ${\mathbf{F}}^{-1}(\mathbf{F}(x))$ 
containing the point $x$. Let $N = \{ L_x \setrule x \in M \} $ and let $\rho : M \rightarrow N: x \mapsto L_x$. 
The integral mapping $\mathbf{F}$ induces the map 
\begin{equation}
\mu : N \rightarrow \mathbf{F}(M) \subseteq {\R }^n: L_x \mapsto \mathbf{F}(x). 
\label{eq-four}
\end{equation}

\noindent \textbf{Claim 5} For every $x \in M$ the connected component $L_x$ of the fiber 
${\mathbf{F}}^{-1}(\mathbf{F}(x))$ containing $x$ is the orbit of the family $X_{\mathbf{F}}$ 
through $x$. \medskip 

\noindent \textbf{Proof.} For each $1 \le i \le n$ let ${\varphi }^{f_i}_t$ be the local flow of the vector field 
$X_{f_i}$ on $(M, \omega )$. Fix $x_0 \in M$. Then $f_j \big( {\varphi }^{f_i}_t(x_0) \big) = f_j(x_0)$ for 
every $1 \le j \le n$. Hence the orbit of $X_{\mathbf{F}}$ through $x_0$ is contained in 
${\mathbf{F}}^{-1}\big( \mathbf{F}(x_0) \big) $. Since orbits of $X_{\mathbf{F}}$ are connected, they 
are the connected components of ${\mathbf{F}}^{-1}\big( \mathbf{F}(x_0) \big) $. \medskip 

To finish the argument we must show that the orbits of $X_{\mathbf{F}}$ are 
open in the fibers of the integral mapping $\mathbf{F}$. Let $O_{x}$ be the orbit of $X_{\boldsymbol{F}}$ through $x$. Suppose that $\mathrm{rank}~\dee \mathbf{F}(x)=k$. By the implicit function theorem,
there is a neighbourhood $U$ of $x$ in $M$ such that $U\cap {\mathbf{F}}^{-1} \big( \mathbf{F}(x) \big)$ is a 
$k$-dimensional submanifold of $M$. On the one hand, since $O_{x}$ is a manifold contained in 
${\mathbf{F}}^{-1}\big( \mathbf{F}(x) \big)$, it follows that its dimension is at most $k$. On the other hand, 
$\mathrm{rank}~\dee  \mathbf{F}(x)=k$ implies that there exist $k$ linear combinations
of vectors $X_{f_{1}}(x),...,X_{f_{n}}(x)$ that are linearly independent at $x$. Therefore, the dimension of the orbit $O_{x}$ is at least $k$. Hence, $\dim O_{x}=k,$ and there exists a neighbourhood $U^{^{\prime }}$ of $x$ in 
$M$ such that $U^{\prime }\subseteq U$ and $U^{\prime }\cap O_{x}=
U^{\prime}\cap {\mathbf{F}}^{-1}\big( \mathbf{F}(x) \big)$ is an open subset of 
${\mathbf{F}}^{-1} \big( \mathbf{F}(x^{0}) \big)$. This holds for every $x\in 
{\mathbf{F}}^{-1}\big( \mathbf{F}(x^{0}) \big) $, which implies that orbits of 
$X_{\mathbf{F}}$ that are contained in ${\mathbf{F}}^{-1}\big( \mathbf{F}(x^{0}) \big)$ are open subsets of 
${\mathbf{F}}^{-1}(\mathbf{F}\big( x^{0})\big) $. \hfill $\square $ \medskip

Claim 5 enables us to identify the space $M/X_{\mathbf{F}}$ of orbits of the family 
$X_{\mathbf{F}}$ of vector fields on $M$ with the space $N$ of connected components of the fibers 
of the integral mapping $\mathbf{F}: M \rightarrow \mathbf{F}(M) \subseteq {\R }^n$. The identification 
$M/X_{\mathbf{F}} =N$ leads to the identification of the projection map $\pi : M \rightarrow M/X_{\mathbf{F}}$ with the 
map $\rho : M \rightarrow N$. In papers on reduction of symmetries in Hamiltonian systems, 
$M/X_{\mathbf{F}}$ is call the orbit space and $\pi :M \rightarrow M/X_{\mathbf{F}}$ the orbit map, 
see \cite{cushman-bates}. In papers on completely integrable Hamiltonian systems, the space 
$N$ of connected components of the fibers of the integral map $\mathbf{F}$ is called the base 
space, see \cite{ratiu-wacheux-zung}. \medskip 

Since $\mathbf{F}$ is constant on the orbits of $X_{\mathbf{F}}$, which are connected components 
of the fibers of $\mathbf{F}$, it follows that the integral map $\mathbf{F}: M \rightarrow \mathbf{F}(M) \subseteq 
{\R }^n$ factors into the composition of $\pi : M \rightarrow M/X_{\mathbf{F}}$ and the map 
\begin{equation}
\mu : M/X_{\mathbf{F}} = N \rightarrow \mathbf{F}(M): L_x \mapsto \mathbf{F}(x). 
\label{eq-five}
\end{equation}
In other words, $\mathbf{F} = \mu \comp \pi $. \medskip 

\noindent \textbf{Corollary 5a} The mapping 
\begin{displaymath}
\mu : \big( M/X_{\mathbf{F}} = N, C^{\infty}(M/X_{\mathbf{F}}) \big)  \rightarrow \big( \mathbf{F}(M), C^{\infty}_i(\mathbf{F}(M)) \big)
\end{displaymath}
is smooth. \medskip 

\noindent \textbf{Proof.} Suppose that $g \in C^{\infty}_i(\mathbf{F}(M))$. From claim 4, we get 
${\mathbf{F}}^{\ast }g \in C^{\infty}(M)$. Clearly ${\mathbf{F}}^{\ast }g = g \comp \mathbf{F}$ is constant on 
the fibers of $\mathbf{F}$. Since the orbits of $X_{\mathbf{F}}$ are connected subsets in the fibers of 
$\mathbf{F}$, it follows that ${\mathbf{F}}^{\ast }g \in {C^{\infty}(M)}^{\mathbf{F}}$. Hence there is a function 
$h \in C^{\infty}(M/X_{\mathbf{F}})$ such that ${\mathbf{F}}^{\ast }g = {\pi }^{\ast }h$. The equality 
$\mathbf{F} = \mu \comp \pi $ yields
\begin{displaymath}
{\pi }^{\ast }h = (\mu \comp \pi )^{\ast }g = {\pi }^{\ast }({\mu }^{\ast }g). 
\end{displaymath}
Since ${\pi }^{\ast }: C^{\infty}(M/X_{\mathbf{F}}) \rightarrow {C^{\infty}(M)}^{\mathbf{F}}$ is bijective, it follows that 
${\mu }^{\ast }g = h \in C^{\infty}(M/X_{\mathbf{F}})$. Hence the mapping $\mu $ is smooth. \hfill $\square $ \medskip 

Suppose that every fiber of the integral map $\mathbf{F}$ of an integrable Hamiltonian system 
$(M, \omega , \mathbf{F})$ is \emph{connected}. Then the smooth mapping 
$\mu : M/X_{\mathbf{F}} \rightarrow \mathbf{F}(M)$ is bijective with inverse 
${\mu }^{-1}: \mathbf{F}(M) \rightarrow M/X_{\mathbf{F}}$. Define a differential structure 
$C^{\infty}\big( \mathbf{F}(M) \big) $ on $\mathbf{F}(M)$ by $C^{\infty}\big( \mathbf{F}(M) \big) = 
({\mu }^{-1})^{\ast }C^{\infty}(M/X_{\mathbf{F}})$. By construction \medskip 

\noindent \textbf{Claim 6} The map 
\begin{equation}
\mu : \big( M/X_{\mathbf{F}}, C^{\infty}(M/X_{\mathbf{F}}) \big) \rightarrow \big( \mathbf{F}(M), C^{\infty}(\mathbf{F}(M)) \big) 
\label{eq-six}
\end{equation}
is a diffeomorphism of differential spaces. \medskip 

The topology $\mathcal{T}$ on $\mathbf{F}(M)$ coming from the differential structure 
$C^{\infty}(\mathbf{F}(M))$ is the same as the topology $\mathcal{S}$ on $M/X_{\mathbf{F}}$ 
coming from the differential structure $C^{\infty}(M/X_{\mathbf{F}})$ because the mapping 
$\mu $ (\ref{eq-six}) is a diffeo and hence homeomorphism. The topology $\mathcal{S}$ is the 
same as the topology ${\mathcal{T}}_1$ on $\mathbf{F}(M)$ coming from the differential structure 
$C^{\infty}_i(\mathbf{F}(M))$, since the mapping $\mu : (M/X_{\mathbf{F}}, \mathcal{S}) \rightarrow 
\big( \mathbf{F}(M), {\mathcal{T}}_1 \big) $ is a continuous bijective map onto a locally compact 
Hausdorff space and thus is a homeomorphism. Consequently, the topologies $\mathcal{T}$ and 
${\mathcal{T}}_1$ on $\mathbf{F}(M)$ are the same. However, the differential spaces 
$\big( \mathbf{F}(M), C^{\infty}(\mathbf{F}(M)) \big) $ and $\big( \mathbf{F}(M), C^{\infty}_i(\mathbf{F}(M)) \big) $ 
are \emph{not} diffeomorphic, since the identity map ${\mathrm{id}}_{\mathbf{F}(M)}$, which is 
the composition of the smooth maps ${\mu }^{-1}$ from $\big( \mathbf{F}(M), C^{\infty}(\mathbf{F}(M)) \big)$ to 
$\big( M/X_{\mathbf{F}}, C^{\infty}(M/X_{\mathbf{F}}) \big) $ and the map $\mu $ from 
$ \big( M/X_{\mathbf{F}}, C^{\infty}(M/X_{\mathbf{F}}) \big) $ to $\big( \mathbf{F}(M), C^{\infty}(\mathbf{F}(M)) \big) $, is a smooth map, whose \linebreak 
inverse is not necessarily smooth as $\mathbf{F}(M)$. \medskip 

\noindent \textbf{Corollary 6a} If $\mathbf{F}(M)$ is a closed subset of ${\R }^n$, then the identity mapping 
\begin{displaymath}
{\mathrm{id}}_{\mathbf{F}(M)}: \big( \mathbf{F}(M), C^{\infty}(\mathbf{F}(M)) \big) 
\rightarrow \big( \mathbf{F}(M), C^{\infty}_i(\mathbf{F}(M)) \big) 
\end{displaymath}
is a diffeomorphism of differential spaces. \medskip 

\noindent \textbf{Proof.} Let $g \in C^{\infty}(\mathbf{F}(M))$. Because $\mathbf{F}(M)$ is a closed 
subset of ${\R }^n$, by the Whitney extension theorem \cite{whitney} there is a smooth function 
$G$ on ${\R }^n$ such that $g = G|_{\mathbf{F}(M)}$. Hence $g \in C^{\infty}_i(\mathbf{F}(M))$. 
So $C^{\infty}(\mathbf{F}(M)) \subseteq C^{\infty}_i(\mathbf{F}(M))$. 
By the above discussion, we have ${\mathrm{id}}^{\ast }_{\mathbf{F}(M)}C^{\infty}_i(\mathbf{F}(M)) \subseteq 
C^{\infty}(\mathbf{F}(M)) \subseteq C^{\infty}_i(\mathbf{F}(M))$. Let $g \in C^{\infty}_i(\mathbf{F}(M))$. Then 
for every $y \in \mathbf{F}(M)$ we have $g(y) = g\big( {\mathrm{id}}_{\mathbf{F}(M)}(y) \big) $ $= 
({\mathrm{id}}^{\ast }_{\mathbf{F}(M)}g)(y)$. So $C^{\infty}_i(\mathbf{F}(M)) \subseteq 
{\mathrm{id}}^{\ast }_{\mathbf{F}(M)}C^{\infty}_i(\mathbf{F}(M))$. Thus 
${\mathrm{id}}^{\ast }_{\mathbf{F}(M)}C^{\infty}_i(\mathbf{F}(M)) $ $= C^{\infty}_i(\mathbf{F}(M))$, which implies 
$C^{\infty}_i(\mathbf{F}(M)) = C^{\infty}(\mathbf{F}(M))$. In other words, the mapping ${\mathrm{id}}_{\mathbf{F}(M)}$ is a diffeomorphism of differential spaces.  \hfill $\square $ \medskip 

\noindent \textbf{Corollary 6b} If $\mathbf{F}(M)$ is a closed subset of ${\R }^n$, then the map 
\begin{displaymath}
\mu : \big( M/X_{\mathbf{F}}, C^{\infty}(M/X_{\mathbf{F}}) \big) \rightarrow \big( \mathbf{F}(M), C^{\infty}_i(\mathbf{F}(M)) \big) 
\end{displaymath} 
is a diffeomorphism.

\end{document}